# POISSON–DIRICHLET DISTRIBUTION FOR RANDOM BELYI SURFACES

By Alex Gamburd[1]

*University of California, Santa Cruz*

*Dedicated to the memory of Bob Brooks*

Brooks and Makover introduced an approach to studying the global geometric quantities (in particular, the first eigenvalue of the Laplacian, injectivity radius and diameter) of a "typical" compact Riemann surface of large genus based on compactifying finite-area Riemann surfaces associated with random cubic graphs; by a theorem of Belyi, these are "dense" in the space of compact Riemann surfaces. The question as to how these surfaces are distributed in the Teichmüller spaces depends on the study of oriented cycles in random cubic graphs with random orientation; Brooks and Makover conjectured that asymptotically normalized cycle lengths follow Poisson–Dirichlet distribution. We present a proof of this conjecture using representation theory of the symmetric group.

**1. Introduction.** The study of the first eigenvalue of the Laplace operator on compact Riemann surfaces of increasing genus has received considerable attention over the last thirty years; see [18] and references therein. On the one hand, we have a celebrated theorem of Selberg [48] (see [36] and [32] for refined estimates toward Selberg's conjecture) and its generalization by Sarnak and Xue [47] (see [29] and [26] for related results), asserting that the first eigenvalue of the congruence surfaces of arbitrary genus is bounded away from zero; on the other hand there are examples due to Selberg [48], Randol [44] and Buser [16], showing that, in general, the first eigenvalue can be made arbitrarily small. The surfaces in the latter examples are "long and thin," so in some sense live on the boundary of Teichmüller spaces, and it is

Received February 2005; revised July 2005.
[1]Supported in part by the NSF postdoctoral fellowship and by the NSF Grant DMS-05-01245.

a fascinating question to determine what happens for a "typical" Riemann surface of large genus.

The idea of using cubic graphs to study the first eigenvalue of Riemann surfaces originated in the work of Buser [15, 17]. As we discuss in Section 3, the behavior of the first eigenvalue of the discrete Laplacian on a random cubic graph is understood rather well. In [14], Brooks and Makover introduced an approach to studying the first eigenvalue of the Laplacian of a "typical" compact Riemann surface of large genus based on compactifying finite-area Riemann surfaces associated with random cubic graphs with random orientation; we review the Brooks–Makover construction in Section 2. To each cubic graph $\Gamma$ with an orientation $\mathcal{O}$, they associate two Riemann surfaces: $S^O(\Gamma, \mathcal{O})$, a finite-area noncompact surface, and $S^C(\Gamma, \mathcal{O})$, a compact surface. The surface $S^O(\Gamma, \mathcal{O})$ is an orbifold cover of $\mathbb{H}/PSL(2,\mathbb{Z})$ described by $(\Gamma, \mathcal{O})$ and therefore shares some of the global geometric properties with the graph $\Gamma$. The compact surface $S^C(\Gamma, \mathcal{O})$ is a conformal compactification of $S^O(\Gamma, \mathcal{O})$; Brooks and Makover proved that almost always the global geometry of $S^C(\Gamma, \mathcal{O})$ is controlled by the geometry of $S^O(\Gamma, \mathcal{O})$. Moreover, according to Belyi's theorem [4], the surfaces $S^C(\Gamma, \mathcal{O})$ are precisely the Riemann surfaces which can be defined over some number field and so form a "dense" set in the space of all Riemann surfaces.

The question as to how these surfaces are distributed in the Teichmüller spaces depends on the study of oriented cycles in random cubic graphs with random orientation; Brooks and Makover conjectured that asymptotically normalized cycle lengths follow Poisson–Dirichlet distribution. We recall the definition of Poisson–Dirichlet distribution [3]. Let $B_1, B_2, \ldots$ be independent random variables uniformly distributed on $[0, 1]$. Define $G = (G_1, G_2, \ldots)$ as follows:

$$G_1 = B_1; \qquad G_2 = (1-B_1)B_2; \qquad G_i = (1-B_1)(1-B_2)\cdots(1-B_{i-1})B_i.$$

The random sequence $G$ can be viewed as a description of a random breaking of a stick of unit length into an infinite sequence of subintervals. A stick of length $B_1$ is broken off at the left, which leaves a piece of length $1 - B_1$. From this, a piece of length $(1-B_1)(1-B_2)$ is broken off, and so on. $G$ is a distribution on the set

$$\Omega = \left[ x \in R^\infty : x_1, x_2, \cdots \geq 0, \sum_{i=1}^\infty x_i = 1 \right].$$

The ranked version of $G$, $(G_{(1)}, G_{(2)}, \ldots)$, where $G_{(1)} \geq G_{(2)} \geq \ldots$, has Poisson–Dirichlet distribution. If, instead of uniform distribution, $B_i$ has beta-$(1, \theta)$ density $\theta(1-x)^{\theta-1}$ on $[0,1]$ with $\theta > 0$, the resulting distribution is called *Poisson–Dirichlet distribution with parameter $\theta$*. Poisson–Dirichlet distribution arises in a great variety of problems; see [3, 42] and references therein. In a recent breakthrough work [21], Diaconis, Mayer-Wolf, Zeitouni and Zerner proved a conjecture of Vershik [52] asserting that Poisson–Dirichlet distribution is the unique invariant distribution for uniform split-merge transformations.

As we discuss in Section 3, the distribution of cycle lengths for random regular graphs in any of the standard models is a rather difficult problem; adding random orientation, at first, seems only to complicate matters (see the example in Section 2). However, as we explain at the end of Section 3, this extra randomness turns out to help us: by giving a permutational model for a random regular graph with random orientation, we can convert the problem into one involving the symmetric group. More precisely, the distribution of oriented cycles in a random $k$-regular graph on $n$ vertices with random orientation is the same as the distribution of cycles in the permutation $\beta\alpha$, where $\beta$ is chosen uniformly on the conjugacy class consisting of the product of $k$-cycles, and $\alpha$ is chosen with uniform probability on the conjugacy class consisting of the product of 2-cycles in the symmetric group $S_N$ with $N = nk$. In Section 4, using the Diaconis–Shahshahani upper bound lemma [22], the estimate on the number of rim hook tableaux by Fomin and Lulov [23] and representation theory of the symmetric group (in particular, the hook length formula and the Murnaghan–Nakayama rule), we show that as $n \to \infty$ the distribution of $\beta\alpha$ converges to uniform distribution. We then invoke what is perhaps the oldest occurrence of Poisson–Dirichlet distribution— the distribution of normalized cycle lengths for a random permutation in $S_n$ as $n$ tends to infinity [49, 57, 58]—to prove the conjecture of Brooks and Makover.

It turns out that the number of oriented cycles in random cubic graphs with random orientation was also studied by Pippenger and Schleich [43] in connection with topological characteristics of random surfaces generated by cubic interactions. The surfaces considered by Pippenger and Schleich are obtained by taking $3n$ arcs of an even number of oriented triangles and randomly identifying them in pairs respecting the orientation; these surfaces arise in various contexts in two-dimensional quantum gravity and as world sheets in string theory. Random cubic graphs with random orientation provide an alternative way of constructing the surfaces in [43]. Denoting by $l(n)$ the number of oriented cycles in a random cubic graph on $n$ vertices with a random orientation, Pippenger and Schleich proved that $\mathbb{E}(l(n)) = \log n + O(1)$ and $\mathrm{Var}(l(n)) = O(\log n)$. Further, based on empirical study of 10,000 random surfaces, each constructed from 80,000 triangles, they conjectured that

(1.1) $$\mathbb{E}(l(n)) = \log(3n) + \gamma + o(1)$$

and

(1.2) $$\mathrm{Var}(l(n)) = \log(3n) + \gamma - \pi^2/6 + o(1),$$

where $\gamma = 0.5772\ldots$ is Euler's constant. It was pointed out by Nicholas Pippenger that (1.1) and (1.2) follow from the main theorem of this paper and a priori bounds on the first few moments of $l(n)$ obtained by the methods of [43]. Another consequence of the main theorem (Corollary 5.2) is that the expected area of the largest embedded ball in a random Belyi surface converges to $\frac{0.62}{2\pi}$ of the total surface area. We conclude Section 5 by briefly discussing the conjectured Tracy–Widom distribution for the second largest eigenvalue of random regular graphs; we hope that our results will be useful in approaching this fascinating open problem.

**2. Belyi surfaces.** In [4], Belyi proved a remarkable result asserting that a Riemann surface $S$ can be defined over the field of algebraic numbers $\overline{\mathbb{Q}}$ if and only if there exists a covering $f : S \to \overline{\mathbb{C}}$ unramified outside $\{0, 1, \infty\}$. We call such surfaces *Belyi surfaces*. In this section we review the Brooks–Makover construction of Belyi surfaces from cubic graphs. We remark that Mulase and Penakava [38] have given a very interesting alternative construction of Belyi surfaces; in their construction, the edges of the graphs are allowed to have variable lengths.

Let $\Gamma$ be a cubic graph. An orientation $\mathcal{O}$ on the graph $\Gamma$ is an assignment for each vertex $v$ of $\Gamma$ of a cyclic ordering of the half-edges incident to that vertex. Given a cubic graph on $n$ vertices, it will have $2^n$ different orientations. A *left-hand-turn path* (LHT path) on $\Gamma$ is a closed path on $\Gamma$ such that, at each vertex, the path turns left in the orientation $\mathcal{O}$.

Given a pair $(\Gamma, \mathcal{O})$, we construct a finite-area Riemann surface $S^O(\Gamma, \mathcal{O})$ as follows. We take the ideal hyperbolic triangle $T$ with vertices $0, 1$ and

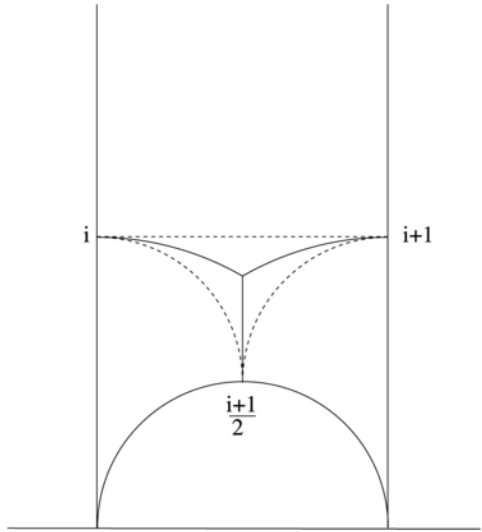

FIG. 1. *The marked ideal triangle $T$.*

∞ shown in Figure 1. The solid lines in Figure 1 are geodesics joining the points $i, i+1$ and $\frac{i+1}{2}$ with the point $\frac{1+i\sqrt{3}}{2}$, while the dotted lines are horocycles joining pairs of points from the set $\{i, i+1, \frac{i+1}{2}\}$. We may think of these points as "midpoints" of the corresponding sides of the ideal triangles, even though the sides are of infinite length. We may also think of the three solid lines as edges of a graph emanating from a vertex. We may then give them the cyclic ordering $(i, i+1, \frac{i+1}{2})$. Given a cubic graph $\Gamma$ with an orientation $\mathcal{O}$, we construct $S^O(\Gamma, \mathcal{O})$ by associating to each vertex an ideal triangle, and gluing neighboring triangles. We glue two copies of $T$ along the corresponding sides, subject to the following two conditions:

(a) the "midpoints" of the two sides are glued together, and
(b) the gluing preserves the orientation of the two copies of $T$.

The conditions (a) and (b) determine the gluing uniquely. It is easily seen that the surface $S^O(\Gamma, \mathcal{O})$ is a complete Riemann surface with a finite area equal to $\pi n$, where $n$ is the number of vertices of $\Gamma$.

It is easy to see that the horocycle pieces on each $T$ glue together to give closed horocycles about a cusp of $S^O(\Gamma, \mathcal{O})$. The length of each closed horocycle is precisely the length of corresponding LHT path. The surface $S^C(\Gamma, \mathcal{O})$ is the conformal compactification of $S^O(\Gamma, \mathcal{O})$. The oriented graph $(\Gamma, \mathcal{O})$ describes $S^O(\Gamma, \mathcal{O})$ as a covering space of $\mathbb{H}^2/PSL_2(\mathbb{Z})$, with each vertex being an orbifold point of order 3.

We are now ready to state the following result, whose proof was sketched by Brooks and Makover in [14]:

LEMMA 2.1.  *$S$ is a Belyi surface if and only if $S = S^C(\Gamma, \mathcal{O})$ for some cubic graph $\Gamma$.*

PROOF. We first show that if $G$ is a torsion-free finite index subgroup of $PSL_2(\mathbb{Z})$, then $\mathbb{H}^2/G = S^C(\Gamma, \mathcal{O})$ for some cubic graph $\Gamma$. Indeed, we can take as a fundamental domain for $PSL_2(\mathbb{Z})$ a set $F = \{0 < \Re(z) < 1, |z| > 1, |z-1| > 1\}$. Three copies of $F$ fit together to give a marked ideal triangle as presented in Figure 1; they are transformed by means of an elliptic element $w$ of order 3 in $PSL_2(\mathbb{Z})$.

Now since $G$ is torsion free, $F$, $w(F)$ and $w^2(F)$ are not equivalent under $G$, and so they can all be included in the fundamental domain of $G$. In particular, there is a fundamental domain of $G$ such that consisting of the copies of ideal triangle $T$. The graph dual to this triangulation, together with the boundary pairings and orientation of $\mathbb{H}^2/G$, is exactly the pair $(\Gamma, \mathcal{O})$.

Now since $S$ is a Belyi surface if and only if one can find finitely many points $\{p_1, \ldots, p_l\}$ on $S$ such that $S - \{p_1, \ldots, p_l\}$ is isomorphic to $\mathbb{H}^2/G$, where $G$ is a torsion-free finite index subgroup of $PSL_2(\mathbb{Z})$ [31], the lemma is proved. □

We define probability on the space of oriented graphs with $n$-vertices $(\Gamma_n, \mathcal{O})$ as follows. We pick a random cubic graph with $n$ vertices using the Bollobas model, described in the next section, then we pick an orientation $\mathcal{O}$ with equal probability from all $2^n$ possible orientations on the given graph. If $Q$ is a property of graphs, we denote by $\text{Prob}_n[Q]$ the probability that an oriented graph $(\Gamma_n, \mathcal{O})$ picked from our probability space has property $Q$. Brooks and Makover proved [14] the following result on Belyi surfaces constructed from random cubic graphs with random orientation:

THEOREM 2.1. *There exist positive constants $C_1$, $C_2$, $C_3$ and $C_4$ such that, as $n \to \infty$:*

(a) *The first eigenvalue $\lambda_1(S^C(\Gamma, \mathcal{O}))$ satisfies*

$$Prob_n[\lambda_1(S^C(\Gamma, \mathcal{O})) \geq C_1] \to 1.$$

(b) *The Cheeger constant $h(S^C(\Gamma, \mathcal{O}))$ satisfies*

$$Prob_n[h(S^C(\Gamma, \mathcal{O})) \geq C_2] \to 1.$$

(c) *The shortest geodesic* $\text{syst}(S^C(\Gamma, \mathcal{O}))$ *satisfies*

$$Prob_n[\text{syst}(S^C(\Gamma, \mathcal{O})) \geq C_3] \to 1.$$

(d) *The diameter* $\text{diam}(S^C(\Gamma, \mathcal{O}))$ *satisfies*

$$Prob_n[\text{diam}(S^C(\Gamma, \mathcal{O})) \leq C_4 \log(\text{genus}(S^C(\Gamma, \mathcal{O})))] \to 1.$$

The theorem is proved, roughly, as follows. First, one establishes that, with probability tending to 1 as $n \to \infty$, a property in question holds for random cubic graphs [see Section 3 for results pertaining to parts (a) and (b)]. Then, using the fact that $(\Gamma, \mathcal{O})$ describes $S^O(\Gamma, \mathcal{O})$ as an orbifold covering, one transfers this information to open surfaces $S^O(\Gamma, \mathcal{O})$, using the results of Brooks in [8, 9]. One then transfers the desired property to the surfaces $S^C(\Gamma, \mathcal{O})$ by using the Ahlfors–Schwarz–Wolpert lemma as developed by Brooks in [10, 11] and extended by Brooks and Makover in [12, 13, 14].

The topology of the surface can be read off from $(\Gamma, \mathcal{O})$, using LHT paths. In particular, the genus is given by Euler's formula; as is customary when using Euler's formula, we refer to an LHT path as a *face* of the oriented graph.

Let $n$ be the number of vertices of $\Gamma$ and let $l(\Gamma, \mathcal{O})$ be the number of disjoint left-hand-turn paths. For a cubic graph, the number of edges is $\mathcal{E}(\Gamma) = \frac{3n}{2}$ and the number of faces is $\mathcal{F} = l(\Gamma, \mathcal{O})$. Therefore, by Euler's formula, we have:

$$(2.1) \qquad \text{genus}(S^C(\Gamma_n, \mathcal{O})) = \text{genus}(S^O(\Gamma_n, \mathcal{O})) = 1 + \frac{n - 2l}{4}.$$

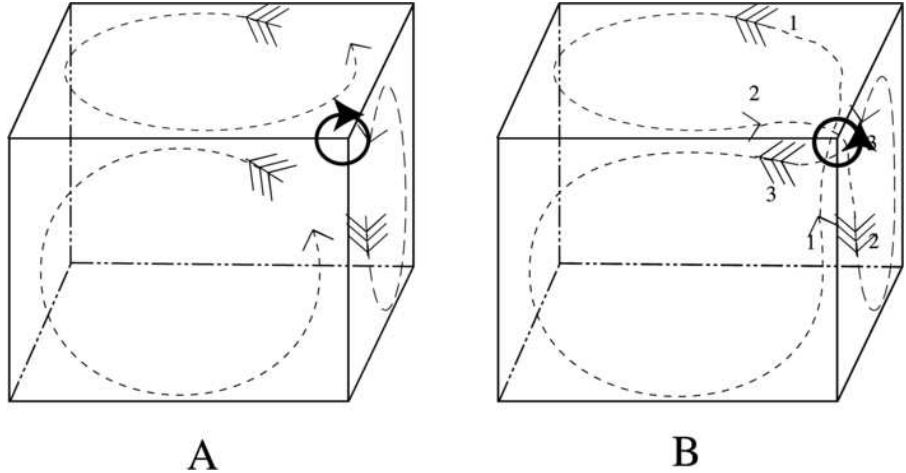

Fig. 2. *Changing orientation on the cube.*

According to Euler's formula, in order to find the genus of a surface constructed from a cubic graph with $n$-vertices, we need to estimate the number of faces (i.e., left-hand-turn paths). As explained in [14], the length of the largest LHT path $L$ also determines the area of the largest embedded ball.

As we discuss in Section 3, computing the distribution of closed paths in random regular graphs is a very difficult problem. An additional complication is that left-hand-turn paths are not necessarily simple closed paths on the graph, but can self-intersect in a very complex pattern. We would like to iliustrate these complications with the simple example of the 1-skeleton of the cube, with the usual orientation (Figure 2A). We have $\mathcal{F} = 6$, and all the faces are simple paths of the graph, therefore the genus is 0. In Figure 2B we changed the orientation of the upper right vertex. With the new orientation, the three simple faces that were previously adjacent to the upper right vertex are now joined to one composite face, while the other three faces are unchanged, hence $\mathcal{F} = 4$, and the surface has genus 1.

**3. Random regular graphs.** In this section we briefly review the pertinent facts on random $k$-regular graphs; see Wormald's paper [59] for an excellent survey. Given a $k$-regular graph $\mathcal{G}$ and a subset $X$ of $V$, the *expansion* of $X$, $c(X)$, is defined to be the ratio $|\partial(X)|/|X|$, where $\partial(X) = \{y \in \mathcal{G} : \text{distance}(y, X) = 1\}$. The *expansion coefficient* of a graph $\mathcal{G}$ is an analogue

of Cheeger's constant for Riemann surfaces and is defined as follows:

$$c(\mathcal{G}) = \inf\{c(X) | |X| < \tfrac{1}{2}|\mathcal{G}|\}. \tag{3.1}$$

A family of $k$-regular graphs $X_{n,k}$ forms a family of $C$-expanders [34, 46] if there is a fixed positive constant $C$ such that

$$\liminf_{n \to \infty} c(X_{n,k}) \geq C. \tag{3.2}$$

The *adjacency matrix* of $\mathcal{G}$, $A(\mathcal{G})$, is the $|\mathcal{G}|$ by $|\mathcal{G}|$ matrix, with rows and columns indexed by vertices of $\mathcal{G}$, such that the $x,y$ entry is 1 if and only if $x$ and $y$ are adjacent, and is 0 otherwise. For a $k$-regular graph, the adjacency matrix is related to the combinatorial Laplacian $\Delta$ by $A = kI - \Delta$. Using the discrete Cheeger–Buser inequality, condition (5.1) can be rewritten in terms of the second largest eigenvalue of the adjacency matrix $A(\mathcal{G})$ as follows:

$$\limsup_{n \to \infty} \lambda_2(A_{n,k}) < k. \tag{3.3}$$

In 1973, Pinsker [41] observed that a random regular graph is a good expander. This corresponds to the following fact about random matrices: a random symmetric matrix of size $N$ with $k$ ones in each row and column and all other entries zero has its biggest eigenvalue equal to $k$, but its next eigenvalue will be bounded away from $k$ by a fixed amount independent of $N$. The result of Pinsker on expansion coefficients of random regular graphs was considerably strengthened by Bollobás [7], who also introduced a widely used configuration model for random regular graphs [6]. In this model, random $k$-regular graphs on $N$ vertices are represented as the images of so-called *configurations*. Let $W = \bigcup_{j=1}^{n} W_j$ be a fixed set of $2m = nd$ vertices, where $|W_j| = d$. A configuration $F$ is a partition of $W$ into $m$ pairs of vertices called *edges* of $F$. Clearly there are

$$N = N(m) = \binom{2m}{2}\binom{2m-2}{2}\cdots\binom{2}{2} \Big/ m! = \frac{(2m)_m}{2^m} \tag{3.4}$$

configurations. [We write $(a)_b = a(a-1)\cdots(a-b+1)$.]

Let $\Phi$ be a set of configurations; we turn it into a probability space by giving all configurations the same probability. We now define a map $\Phi \to \mathcal{G}_{n,k}$ as follows. Given a configuration $F$, let $\phi(F)$ be the graph with vertex set $V = 1, \ldots, n$ in which $ij$ is an edge iff $F$ has a pair with one end in $W_i$ and the other in $W_j$. Every $G \in G_{n,d}$ is the image of $G = \phi(F)$ for $(d!)^n$ configurations. The number of configurations containing a given fixed set of $l$ edges is

$$\begin{aligned}
N_l(m) &= \binom{2m-2l}{2}\binom{2m-2l-2}{2}\cdots\binom{2}{2} \Big/ (m-l)! \\
&= \frac{N(m)}{(2m-1)(2m-3)\cdots(2m-2l+1)}.
\end{aligned} \tag{3.5}$$

Using the configuration model, and in particular (3.5), Bollobás proved the following result:

THEOREM 3.1 ([5]). *Let $X_i$ denote the number of closed walks of length $i$ in $\mathcal{G} \in \mathcal{G}_{n,k}$. Then, for $i = o(\log_{k-1} n)$ as $n \to \infty$, random variables $X_i$ converge to independent Poisson random variables with mean $\frac{(k-1)^i}{2i}$.*

Counting cycles of length greater than $\log n$ is substantially more difficult. In a recent breakthrough work [25], following his earlier work in [24], Friedman estimates the number of cycles of length $O(\log^2 n)$ and uses this estimate (among other things) to prove that $k$-regular graphs on $n$ vertices $\mathcal{G}_{n,k}$ are asymptotically Ramanujan: for $k$ fixed and $\varepsilon > 0$, the probability that $\lambda_1(X_{n,k}) \leq 2\sqrt{k-1} + \varepsilon$ tends to 1 as $n \to \infty$.

The bound of $2\sqrt{k-1}$ is optimal in view of the result of Alon–Boppana [1, 35]. We also mention an early result of McKay [37], who showed that spectral density of random $k$-regular graphs converges to Kesten's measure, that is, a measure supported on $[-2\sqrt{k-1}, 2\sqrt{k-1}]$ and given by

$$\nu_k = \frac{k}{2\pi} \frac{\sqrt{4(k-1) - t^2}}{k^2 - t^2}. \tag{3.6}$$

We now introduce the permutational model for random regular graphs with random orientation. Consider a neighborhood of a vertex of a cubic graph. It contains three half-edges incident to it. Denote the set of all half-edges by $H$; $|H| = 2|E| = 3|V|$. Now the cyclic ordering of half-edges at each vertex is specified by a 3-cycle, and a collection of all cyclic orders on the half-edges yields a permutation of $H$ consisting of the product of 3-cycles. The structure of the underlying graph is given by the way half-edges couple to each other to form an edge; this is described by a permutation all of whose cycles are of order 2.

In the Figure 3 we have

$$\beta = (1,3,5)(2,12,8)(4,7,9)(6,10,11),$$
$$\alpha = (1,2)(3,4)(5,6)(7,8)(9,10)(11,12)$$

and

$$\varphi = \beta\alpha = (1,12,6)(2,3,7)(4,5,10)(8,9,11).$$

Let $\mathcal{G} \in \mathcal{G}_{n,k}$ be a $k$-regular graph on $n$ vertices; let $\mathcal{O}$ be its orientation, that is, a cyclic ordering of incident half-edges for each vertex. Associated with a $(\mathcal{G}, \mathcal{O})$, we have a pair of permutations $(\beta, \alpha)$ in $S_N$, where $N = kn$, with cycles of $\beta$ encoding the information about the vertices with their cyclic orientation and $\alpha$ encoding the information about the edges. The cycles of $\varphi = \beta\alpha$ encode the information about the faces.

Now the problem of the distribution of the lengths of LHT paths can be given the following, equivalent, formulation. Denote by $C_r$ the conjugacy class of $A_N$ consisting of the product of $N/r$ disjoint $r$-cycles. Choose $\beta$ with uniform probability on $C_k$ and choose $\alpha$ with uniform probability on $C_2$. Then the distribution of LHT paths in a random $k$-regular graph on $n$ vertices with random orientation is the same as the distribution of cycles in the permutation $\beta\alpha$.

**4. Proof of Brooks–Makover conjecture.** This section is devoted to the proof of the following result:

THEOREM 4.1. *Let $N = kn$ and let $P_r$ denote the probability measure on $A_N$ supported on $C_r$; let $U$ denote the uniform distribution on $A_N$. Then, for $k \geq 3$,*

$$\|P_k * P_2 - U\| \xrightarrow[n \to \infty]{} 0. \tag{4.1}$$

*Here*

$$\|f - g\| = \max_{A \subseteq G} |f(A) - g(A)| \tag{4.2}$$

*is a total variation distance.*

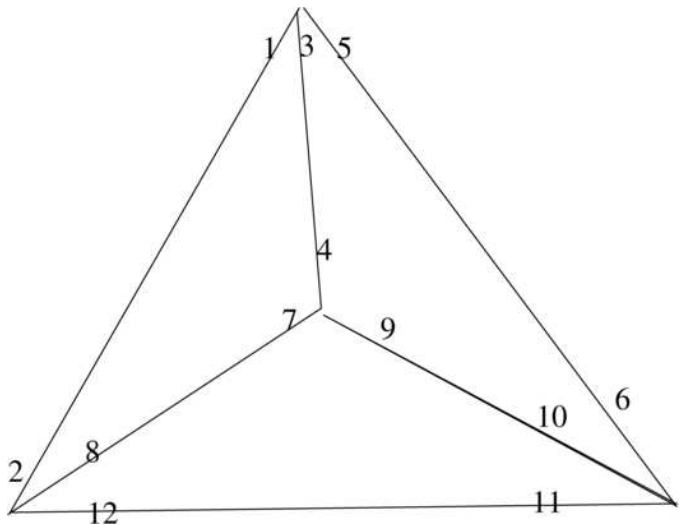

FIG. 3. *Permutational representation.*

COROLLARY 4.1. *The distribution of LHT paths for random regular graphs with random orientation converges to Poisson–Dirichlet distribution; in particular, the Brooks–Makover conjecture is true.*

PROOF. Let $\pi$ denote the cycle type of a permutation; denote by $C^S(\pi)$ the number of permutations with a given cycle type $\pi$ in the symmetric group $S_N$ and by $C^A(\pi)$ the number of permutations with cycle type $\pi$ in the alternating group $A_N$. As is well known,

$$\begin{aligned}C^S(\pi) &= C^S(a_1, \ldots, a_N) \\ &= \frac{N!}{\prod_i i^{a_i(\pi)} a_i(\pi)!}\end{aligned} \qquad (4.3)$$

for $\pi = 1^{a_1} \ldots N^{a_N}$. The cycle index polynomial of $S_N$ is given by

$$p_N(x_1, \ldots, x_N) = \sum_{1a_1 + \cdots + Na_N = N} C^S(a_1, \ldots, a_N) x^{a_1} \cdots x^{a_N}.$$

It is easy to see that the cycle index polynomial of $A_N$, $q_N(x_1, \ldots, x_N)$, is given by

$$q_N(x_1, \ldots, x_N) = p_N(x_1, \ldots, x_N) + p_N(x_1, -x_2, \ldots, (-1)^{N-1} x_N).$$

It follows that

$$C^A(\pi) = \begin{cases} 2C^S(\pi), & \text{if } a_2 + \cdots + a_N \text{ is even,} \\ 0, & \text{otherwise.} \end{cases} \qquad (4.4)$$

Now the condition that $a_2 + \cdots + a_N$ is even, is exactly the condition for a permutation to be in $A_N$ and, subject to this condition, we see that the distribution of cycles is given by exactly the same formula as in the case of the symmetric group. Consequently, the proofs of the asymptotic Poisson–Dirichlet distribution for normalized cycle lengths in the symmetric group apply also in the case of alternating group. The proof of the corollary is completed by applying Theorem 4.1 and the triangle inequality. □

We now turn to the proof of Theorem 4.1. The basic tool is the following result, known as the Diaconis–Shahshahani upper bound lemma; see [19, 20] for a survey of applications and ramifications:

PROPOSITION 4.1 ([22]). *Let $G$ be a finite group and denote by $\hat{G}$ the set of irreducible unitary representations of $G$. Let $P$ be a probability measure on $G$ and denote by $\hat{P}(\rho)$ its Fourier transform at the representation $\rho \in \hat{G}$. Let $U$ be uniform probability measure on $G$. Then*

$$\|P - U\|^2 \leq \tfrac{1}{4} \sum_{\substack{\rho \in \hat{G} \\ \rho \neq \mathrm{id}}} \dim(\rho) \mathrm{tr}(\hat{P}(\rho) \overline{\hat{P}(\rho)}). \qquad (4.5)$$

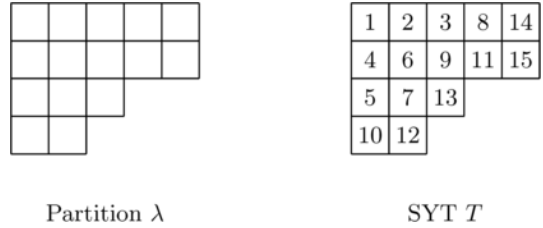

Partition $\lambda$      SYT $T$

Fig. 4.

Here $\|P - U\|$ is the total variation distance defined in (4.2) and 'id' denotes the trivial representation.

We apply Proposition 4.1 to the situation where $G = A_N$ and $P = P_k * P_2$. Since $P_k$ is a class function, $\hat{P}_k(\rho)$ is a scalar matrix given by $\chi^\rho(C_k) I_{\dim(\rho)}$, where $\chi^\rho(C_k)$ is a character of the alternating group associated with representation $\rho$ evaluated on the conjugacy class $C_k$ and $I_{\dim(\rho)}$ is the identity matrix of dimension $\dim(\rho)$. Now the Fourier transform maps the convolution of functions $P_k * P_2$ into their product $\widehat{P_k * P_2} = \widehat{P_k}\widehat{P_2}$ and, consequently, (4.5) implies

$$(4.6) \qquad \|P_k * P_2 - U\|^2 \leq \frac{1}{4} \sum_{\substack{\rho \in \widehat{A_N} \\ \rho \neq \mathrm{id}}} \left( \frac{\chi^\rho(C_k)\chi^\rho(C_2)}{\dim(\rho)} \right)^2.$$

The representation theory of the alternating group $A_N$ is closely allied with the representation theory of the symmetric group $S_N$ [30]. Representations of the symmetric group $S_N$ are labeled by partitions $\lambda \vdash N$. A *partition* $\lambda$ of a nonnegative integer $N$ is a sequence $(\lambda_1, \ldots, \lambda_r) \in \mathbb{N}^r$ satisfying $\lambda_1 \geq \cdots \geq \lambda_r$ and $\sum \lambda_i = N$. We call $|\lambda| = \sum \lambda_i$ the *size* of $\lambda$. The number of parts of $\lambda$ is the length of $\lambda$, denoted $l(\lambda)$. We write $m_i = m_i(\lambda)$ for the number of parts of $\lambda$ that are equal to $i$, so we have $\lambda = \langle 1^{m_1} 2^{m_2} \ldots \rangle$.

The *Young diagram* of a partition $\lambda$ is defined as the set of points $(i,j) \in \mathbb{Z}^2$ such that $1 \leq i \leq \lambda_j$; it is often convenient to replace the set of points above by squares. The *conjugate partition* $\lambda'$ of $\lambda$ is defined by the condition that the Young diagram of $\lambda'$ is the transpose of the Young diagram of $\lambda$; equivalently $m_i(\lambda') = \lambda_i - \lambda_{i+1}$.

A *standard Young tableau* (SYT) of shape $\lambda$ is a filling of the boxes of $\lambda$ with positive integers $1, \ldots, N$ such that the rows and the columns are strictly increasing.

In Figure 4 we exhibit a partition $\lambda = (5,5,3,2) = \langle 1^0 2^1 3^1 5^2 \rangle$ and an SYT $T$ of shape $\lambda$.

We denote by $f^\lambda$ the dimension of representation of $S_N$ indexed by $\lambda$; it is equal to the number of standard Young tableaux of shape $\lambda$ [30]. We denote by $\chi^\lambda(\mu)$ the value of the character indexed by $\lambda$ on the conjugacy class $C_\mu$. For example, for the conjugacy class $C_k$, the corresponding partition is $\mu = k^n$.

The partition $\lambda$ is called *self-associated* if its Young diagram is symmetric with respect to the main diagonal, equivalently if $\lambda = \lambda'$. Now for an irreducible representation of $S_N$ indexed by $\lambda$ which is not self-associated, the restriction of $\lambda$ to $A_N$ is irreducible ([30], page 67). If $\lambda$ *is* self-associated then restriction to $A_N$ splits into two irreducible representations, $\lambda^+$ and $\lambda^-$. The character values of the irreducible representations $\lambda^+$ and $\lambda^-$ can be expressed in terms of character values of $\chi^\lambda(\mu)$; for the conjugacy class $\mu$ not equal to the set of main diagonal hooks of $\lambda$ (see the definition of hook below), it is given by $\frac{\chi^\lambda(\mu)}{2}$. All irreducible representations of $A_N$ are obtained from the irreducible representations of $S_N$ in this way. Consequently, we obtain the following lemma:

LEMMA 4.1. *Notation being as above, we have*

$$\text{(4.7)} \qquad \sum_{\substack{\rho \in \widehat{A_N} \\ \rho \neq \text{id}}} \left( \frac{\chi^\rho(C_k)\chi^\rho(C_2)}{\dim(\rho)} \right)^2 \leq 2 \sum_{\substack{\lambda \in \widehat{S_N} \\ \lambda \neq \langle N \rangle, \langle 1^N \rangle}} \left( \frac{\chi^\lambda(C_k)\chi^\lambda(C_2)}{f^\lambda} \right)^2.$$

PROOF. Indeed, if $\lambda$ is not a self-associated partition, we have $\dim_{A_N}(\lambda) = f^\lambda$ and $\chi^\lambda_{A_N}(\mu) = \chi^\lambda(\mu)$. If $\lambda$ *is* a self-associated partition, then, since conjugacy class $C_k$ is not an array of main diagonal hooks of $\lambda$, we have $\dim_{A_N}(\lambda^\pm) = \frac{f^\lambda}{2}$ and $\chi^{\lambda^\pm}_{A_N}(C_k) = \frac{\chi^\lambda(C_k)}{2}$. This completes the proof of Lemma 4.1. □

We thus have to analyze the sum

$$\text{(4.8)} \qquad \sum_{\substack{\lambda \vdash N \\ \lambda \neq \langle N \rangle, \langle 1^N \rangle}} \left( \frac{\chi^\lambda(C_k)\chi^\lambda(C_2)}{f^\lambda} \right)^2.$$

We briefly recall the basic facts pertaining to $f^\lambda$ and $\chi^\lambda(\mu)$, referring to [30] for more details. Given a diagram $\lambda$ and a square $u = (i, j) \in \lambda$, a *hook* with vertex $u$ is a set of squares in $\lambda$ directly to the right or directly below $u$. We define *hook length* (also referred to as *hook number*) $h(u)$ of $\lambda$ at $u$ by

$$h(u) = \lambda_i + \lambda'_j - i - j - 1.$$

Equivalently, $h(u)$ is the number of squares directly to the right or directly below $u$, counting $u$ itself once. For instance, in Figure 5 we display hook lengths for the partition $\lambda = (5, 5, 3, 2)$.

THEOREM 4.2 (Hook length formula). *Notation being as above, we have*

$$f^\lambda = \frac{N!}{\prod_{u \in \lambda} h(u)}. \tag{4.9}$$

Projecting a hook with vertex $u$ onto the boundary (rim) of $\lambda$ yields a *rim hook*. By definition, a rim hook is a connected skew shape with no $2 \times 2$ square. The *height* of a rim hook, $\operatorname{ht}(R)$, is defined to be one less than its number of rows.

A rim hook tableau $T$ of shape $\lambda$ and type $\mu$, where $\mu = (\mu_1, \ldots, \mu_n)$, is an assignment of positive integers to the squares of $\lambda$ such that every row and column is weakly increasing; the integer $i$ appears $\mu_i$ times, and the set of squares occupied by $i$ forms a rim hook. The *height* of a rim hook tableau is defined to be the sum of the heights of rim hooks appearing in $T$.

In Figure 6, we exhibit a rim hook tableaux of shape $\lambda = (5, 5, 3, 2)$ and type $\mu = (4, 4, 4, 3)$.

In the particular case when all the parts of $\mu$ have size $k$, the tableaux described above is referred to as a *$k$-rim hook tableaux*; we denote the number of $k$-rim hook tableaux by $f_k^\lambda$. In particular, $f_1^\lambda = f^\lambda$ is the number of a standard Young tableaux.

THEOREM 4.3 (Murnaghan–Nakayma rule). *Notation being as above, we have*

$$\chi^\lambda(\mu) = \sum_T (-1)^{\operatorname{ht}(T)}, \tag{4.10}$$

*where the sum is over all rim hook tableaux of shape $\lambda$ and type $\mu$.*

| 8 | 7 | 5 | 3 | 2 |
|---|---|---|---|---|
| 7 | 6 | 4 | 2 | 1 |
| 4 | 3 | 1 | | |
| 2 | 1 | | | |

FIG. 5.

| 1 | 1 | 3 | 3 | 4 |
|---|---|---|---|---|
| 1 | 2 | 3 | 4 | 4 |
| 1 | 2 | 3 | | |
| 2 | 2 | | | |

FIG. 6.

For $\mu = k^n$, all signs (4.10) are the same [30], page 82; consequently,

$$|\chi^\lambda(C_k)| = f_k^\lambda. \tag{4.11}$$

In the proof of Theorem 4.1 we will make crucial use of the following estimate for $f_k^\lambda$:

THEOREM 4.4 ([23]). *Suppose $N = kn$. Then*

$$f_k^\lambda \leq \frac{n! k^n}{(kn)!^{1/k}} (f^\lambda)^{1/k} = O(N^{1/2 - 1/(2k)})(f^\lambda)^{1/k}, \tag{4.12}$$

*where the implied constant depends only on $k$ and not on $\lambda$.*

Another crucial ingredient in the proof of Theorem 4.1 is furnished by the following proposition:

PROPOSITION 4.2. *For any $t > 0$, we have*

$$\sum_{\substack{\lambda \vdash N \\ \lambda_1, \lambda_1' \leq N - m}} (f^\lambda)^{-t} = O(N^{-mt}), \tag{4.13}$$

*where the implied constant depends only on $m$.*

We remark that the sums of the form $\sum_{\lambda \vdash N} (f^\lambda)^\beta$ for $\beta > 0$ have been studied by Regev [45] and Vershik [56]. In particular, Regev relates the asymptotic computations of such sums to the matrix integral of the form

$$\int_{-\infty}^{\infty} \cdots \int_{-\infty}^{\infty} \prod_{i<j} |x_i - x_j|^\beta e^{-(\beta/2)(x_1^2 + \cdots + x_N^2)} \, dx_1 \cdots dx_N,$$

this being one of the first hints of the deep connection between random matrices and random permutations; see [40] and references therein for a recent survey.

PROOF OF PROPOSITION 4.2. First, we observe that since $f^\lambda = f^{\lambda'}$, it suffices to prove Proposition 4.2 for the sum

$$\sum_{\substack{\lambda \vdash N \\ \lambda_1' < \lambda_1 \leq N - m}} (f^\lambda)^{-t}.$$

Now we split this sum into three parts:

$$\Sigma_1 = \sum_{\substack{\lambda \vdash N \\ \lambda_1' \leq \lambda_1 \\ 3N/4 < \lambda_1 \leq N - m}} (f^\lambda)^{-t}, \tag{4.14}$$

$$\Sigma_2 = \sum_{\substack{\lambda \vdash N \\ \lambda'_1 \leq \lambda_1 \\ N/8 < \lambda_1 \leq 3N/4}} (f^\lambda)^{-t} \tag{4.15}$$

and

$$\Sigma_3 = \sum_{\substack{\lambda \vdash N \\ \lambda'_1 \leq \lambda_1 \\ \lambda_1 \leq N/8}} (f^\lambda)^{-t}. \tag{4.16}$$

To analyze $\Sigma_1$, we first note that

$$f^\lambda \geq \binom{\lambda_1}{N-\lambda_1} \quad \text{if } \lambda_1 > \frac{N}{2}. \tag{4.17}$$

Indeed, since the product of hook numbers for a partition of size $l$ is at most $l!$, the product of hook numbers in the rows below the top one is at most $(N-\lambda_1)!$. The set of hook numbers in the top row consists of distinct numbers not exceeding $N$. Since the length of the second row is at most $r = N - \lambda_1$, the hook numbers in the rightmost $\lambda_1 - r$ cells of the first row are $1, \ldots, 2\lambda_1 - N$. The product of the remaining $r$ hook numbers in the first row is no greater than $N(N-1)\ldots(N-r+1) = \frac{N!}{\lambda_1!}$. Applying the hook length formula (4.9) then completes the proof of the estimate (4.17).

Continuing to denote $r = N - \lambda_1$, we can now estimate $\Sigma_1$ using (4.17):

$$\Sigma_1 \leq \sum_{m \leq r \leq N/4} \frac{p(r)}{\binom{N-r}{r}^t}, \tag{4.18}$$

where $p(r)$ is the number of partitions of $r$. Since the number of partitions $p(r)$ satisfies the following inequality [2] valid for all $r \geq 1$:

$$p(r) \leq \exp(\pi\sqrt{2r/3}), \tag{4.19}$$

we have

$$\Sigma_1 \leq \sum_{m \leq r \leq N/4} \frac{c_1^{\sqrt{r}}}{\binom{N-r}{r}^t} \tag{4.20}$$

for absolute constant $c_1 = e^{\pi\sqrt{2/3}}$. Let $a_r = \frac{1}{\binom{N-r}{r}}$. We have

$$a_{r+1} = \frac{(N-r)(r+1)}{(N-2r)(N-2r-1)} a_r.$$

For large $N$, as $r$ increases from $m$ to $\frac{N}{4}$ the function $\frac{(N-r)(r+1)}{(N-2r)(N-2r-1)}$ monotonically decreases from roughly $\frac{1}{N}$ to $\frac{3}{4}$. Consequently,

$$\Sigma_1 \leq \left(\frac{1}{(N-m)\ldots(N-2m)}\right)^t \sum_{m \leq r \leq N/4} \frac{c_1^{\sqrt{r}}}{c_2^r},$$

where $c_2 = (\frac{4}{3})^t > 1$. Now

$$\sum_{m \leq r \leq N/4} \frac{c_1^{\sqrt{r}}}{c_2^r} < \sum_{m \leq r \leq \infty} \frac{c_1^{\sqrt{r}}}{c_2^r},$$

which converges by the Cauchy criterion and, consequently, we obtain the desired estimate $\Sigma_1 = O(N^{-mt})$.

To analyze $\Sigma_2$, we first note that

(4.21) $$f^\lambda \geq \frac{(17N/16 - \lambda_1)!}{(N-\lambda_1)!(N/16+16)!} \quad \text{if } \lambda_1 \geq \frac{N}{8}.$$

The proof of this estimate is similar to the proof of the estimate (4.17). The product of hook numbers outside the first row is at most $(N-\lambda_1)!$. Suppose there is a total of $v$ rows. Then

$$\lambda_2 + \cdots + \lambda_v < \frac{7N}{8}.$$

Consequently,

$$\left(\lambda_2 - \frac{N}{16}\right) + \cdots + \left(\lambda_{15} - \frac{N}{16}\right) < 0,$$

and therefore the rows below the 16th one are of size at most $\frac{N}{16}$. That means that the product of hook numbers in the $\frac{N}{16}$ leftmost boxes of the first row is no greater than $(\frac{N}{16}+16)!$. The remaining hook lengths in the first row constitute $\lambda_1 - \frac{N}{16}$ distinct numbers less than or equal to $N$; their product therefore does not exceed $\frac{N!}{(N+N/16-\lambda_1)!}$. Application of the hook length formula (4.9) then completes the proof of the estimate (4.21).

Now as $\lambda_1$ increases from $\frac{N}{4}$ to $\frac{3N}{4}$, the expression on the right-hand side of (4.21) decreases; therefore, in the whole range of the sum $\Sigma_2$ we have

$$f^\lambda \geq \frac{(17N/16 - 3N/4)!}{(N-3N/4)!(N/16+16)!} > \frac{(5N/16)!}{(N/4)!(N/16)!},$$

and the latter expression is greater than $5^{N/16}$ for sufficiently large $N$. Consequently,

$$\Sigma_2 \leq p(N) 5^{-tN/16} \leq c_1^{\sqrt{N}} c_3^{-N} = O(c_4^{-N})$$

for $c_3 = 5^{t/16}$ and some $c_4 > 1$.

Turning now to the sum $\Sigma_3$, we note that

$$(4.22) \qquad f^\lambda \geq \left(\frac{4}{e}\right)^N \qquad \text{if } \lambda_1' \leq \lambda_1 < \frac{N}{8}.$$

Indeed, for any box $u$ in $\lambda$ we have $h(u) \leq \lambda_1 + \lambda_1' < \frac{N}{4}$. Consequently, using the hook length formula (4.9) we obtain

$$f^\lambda > \frac{N!}{(N/4)^N} > \frac{(N/e)^N}{(N/4)^N} > \left(\frac{4}{e}\right)^N,$$

proving the estimate (4.22). Consequently,

$$\Sigma_3 \leq c_1^{\sqrt{N}} c_5^{-N} = O(c_6^{-N})$$

for $c_5 = (\frac{4}{e})^t$ and some $c_6 > 1$.

This completes the proof of Proposition 4.2.  □

Now Theorem 4.1 is proved by combining inequality (4.6) with Lemma 4.1, Theorem 4.4 and Proposition 4.2 with $m = 4$, together with the values of the remaining characters and dimensions, which are computed using Theorem 4.2 and Theorem 4.3 and summarized in Table 1. We denote by $\Lambda$ the set of partitions in the first column of Table 1.

$$\|P_k * P_2 - U\|^2 \leq \frac{1}{4} \sum_{\substack{\rho \in \widehat{A_N} \\ \rho \neq \mathrm{id}}} \left(\frac{\chi^\rho(C_k)\chi^\rho(C_2)}{\dim(\rho)}\right)^2 \qquad \text{by inequality (4.6)}$$

$$\leq \frac{1}{2} \sum_{\substack{\lambda \in \widehat{S_N} \\ \lambda \neq \langle N \rangle, \langle 1^N \rangle}} \left(\frac{\chi^\lambda(C_k)\chi^\lambda(C_2)}{f^\lambda}\right)^2 \qquad \text{by Lemma 4.1}$$

TABLE 1
*Character and dimension values*

| $\lambda$ | $f^\lambda$ | $\|\chi^\lambda(C_2)\|$ | $\|\chi^\lambda(C_3)\|$ | $\|\chi^\lambda(C_4)\|$ | $\|\chi^\lambda(C_k)\|\ k \geq 5$ |
|---|---|---|---|---|---|
| $(N-1,1)$ | $N-1$ | 1 | 1 | 1 | 1 |
| $(N-2,2)$ | $\frac{N(N-3)}{2}$ | $\frac{N}{2}$ | 0 | 1 | 1 |
| $(N-2,1,1)$ | $\frac{(N-1)(N-2)}{2}$ | $\frac{N}{2}+1$ | 1 | 1 | 1 |
| $(N-3,2,1)$ | $\frac{N(N-2)(N-4)}{3}$ | 0 | $\frac{N}{3}+1$ | 0 | 1 |
| $(N-3,1,1,1)$ | $\frac{(N-1)(N-2)(N-3)}{3}$ | $\frac{N}{2}+1$ | $\frac{N}{3}-1$ | 1 | 1 |
| $(N-3,3)$ | $\frac{(N)(N-1)(N-5)}{6}$ | $\frac{N}{2}+2$ | $\frac{N}{3}+1$ | 0 | $\begin{cases} 0, & \text{if } k=5; \\ 1, & \text{if } k>5 \end{cases}$ |

$$= \frac{1}{2} \sum_{\substack{\lambda \vdash N \\ \lambda_1, \lambda'_1 \leq N-4}} \left(\frac{\chi^\lambda(C_k)\chi^\lambda(C_2)}{f^\lambda}\right)^2 + \sum_{\substack{\lambda \vdash N \\ \lambda \in \Lambda}} \left(\frac{\chi^\lambda(C_k)\chi^\lambda(C_2)}{f^\lambda}\right)^2$$

$$= O(N^{3/2-1/k}) \sum_{\substack{\lambda \vdash N \\ \lambda_1, \lambda'_1 \leq N-4}} (f^\lambda)^{-(1-2/k)} + O(N^{-2})$$

by Theorem 4.4 and Table 1

$$= O(N^{-5/2+7/k}) + O(N^{-2}) \qquad \text{by Proposition 4.2;}$$

the first term in the last expression is at most $O(N^{-1/6})$, since $k \geq 3$. This completes the proof of Theorem 4.1.

## 5. Concluding remarks.

5.1. Extensive information is available on the distribution of cycles in random permutations and on Poisson–Dirichlet distribution [3, 42] and we will pursue the exhaustive exploitation of the consequences of Theorem 4.1 in a subsequent publication. Here we note just two immediate corollaries which substantially refine results in [14, 28, 43].

COROLLARY 5.1. *Let $l(n)$ denote the number of LHT paths in a random cubic graph on $n$ vertices with random orientation. Then, as $n \to \infty$,*

(5.1) $$\mathbb{E}(l(n)) = \log(3n) + \gamma + o(1),$$

(5.2) $$\mathrm{Var}(l(n)) = \log(3n) + \gamma - \pi^2/6 + o(1),$$

*where $\gamma = 0.5772\ldots$ is Euler's constant. Further, $\frac{l(n) - \log n}{\sqrt{\log n}}$ converges to standard normal distribution $\mathcal{N}(0,1)$. From (2.1) we obtain that the genus is distributed as $1 + \frac{n}{4} - \mathcal{N}(\log n, \sqrt{\log n})$.*

COROLLARY 5.2. *Let $L(n)$ be the length of the largest LHT path. Then*

(5.3) $$\lim_{n \to \infty} \mathbb{E}\left(\frac{L(n)}{n}\right) = \int_0^\infty \exp\left(-x - \int_x^\infty (e^{-y}/y)\,dy\right) dx \sim 0.6243.$$

*Consequently, the expected area of the largest embedded ball in $S^C(\Gamma, \mathcal{O})$ converges to $\frac{0.62}{2\pi}$ of the total surface area.*

The expression on the right-hand side of (5.3) is due to Shepp and Lloyd [49], who also computed the limiting distribution of $L$.

5.2. Recent numerical experiments of Novikoff [39] present convincing evidence in favor of the following conjecture (see [27] for a discussion of related conjectures and numerical results):

CONJECTURE 5.1 ([39]). The distribution of the second largest eigenvalue of the adjacency matrix of a random regular graph, suitably rescaled, follows Tracy–Widom GOE distribution.

Tracy and Widom [53, 54] computed the limiting distribution function for the largest eigenvalue in the classical Gaussian ensembles; these distribution functions are expressed in terms of a certain Painlevé II function and are now believed to describe new universal limit laws for a wide variety of processes arising in mathematical physics and interacting particle systems [55].

A dramatic consequence of Conjecture 5.1 would be that the probability of a random regular graph being Ramanujan approaches 0.52 as the size of the graph tends to infinity, corresponding to the skewness in the Tracy–Widom GOE distribution.

To approach Conjecture 5.1 following the method Sinai and Soshnikov [50] and Soshnikov [51] in his breakthrough proof of the universality at the edge of the spectrum in Wigner matrices, one needs precise information for the number of closed walks of size up to $n^{2/3}$, where $n$ is the size of the graph. We hope that the results of this paper will be useful in pursuing such an approach.

**Acknowledgments.** I started working on this problem while visiting Bob Brooks at the Technion in August 2001; I would like to thank the Department of Mathematics at Technion for its hospitality. The interplay between spectral geometry of surfaces and graphs, particularly the issue of the first eigenvalue, has been a recurrent theme in Bob's work; his papers, written with deep sympathy for the reader, and his continuing support and encouragement have been very important to me. The results presented in this paper were obtained in the spring of 2002 and I communicated them to Bob a few months before his untimely death in September 2002. I lectured on them in June 2002 at the spectral geometry conference in Lexington, Kentucky and at the Brooks Memorial Conference in Technion in December 2003; I would like to thank the organizers of these conferences. After my talk at the Brooks Memorial Conference, Alex Lubotzky told me that Martin Liebeck and Aner Shalev independently obtained a result similar to Proposition 4.2; I would like to thank Aner Shalev for sending me their preprint [33]. I am grateful to Persi Diaconis, Eran Makover, Nicholas Pippenger, Anatoly Vershik, Ofer Zeitouni, Bálint Virág and Martin Zerner for interest in this work and for stimulating discussions.

DEPARTMENT OF MATHEMATICS
UNIVERSITY OF CALIFORNIA
1156 HIGH STREET
SANTA CRUZ, CALIFORNIA 95064
USA
E-MAIL: agamburd@ucsc.edu